\documentclass[11pt,twoside]{article}
\usepackage{latexsym}
\usepackage{amssymb,amsbsy,amsmath,amsfonts,amssymb,amscd}
\usepackage{epsfig, graphicx}
\usepackage{color}
\usepackage[pageanchor=true]{hyperref}
\usepackage{pdfsync}
\usepackage[active]{srcltx}
\usepackage[latin1]{inputenc}
\usepackage[francais]{babel}
\newcommand{\C}{\mathbb{C}}
\newcommand{\R}{\mathbb{R}}

\def\Box{\hfill\rule{2.5mm}{2.5mm}}

\newcommand {\f}   {\frac}

\newcommand{\beq}{\begin{equation}}
\newcommand{\beqa}{\begin{eqnarray}}
\newcommand{\bea} {\begin{array}{ll}}
\newcommand{\beqan}{\begin{eqnarray*}}
\newcommand{\eeq}{\end{equation}}
\newcommand{\eeqa}{\end{eqnarray}}
\newcommand{\eeqan}{\end{eqnarray*}}
\newcommand{\eea} {\end{array}}
\newcommand{\boxs}{\hfill\rule{2.5mm}{2.5mm}}
\newtheorem{cor}{Corollary}[section]

\newtheorem{lem}[cor]{Lemma}
\newtheorem{prop}[cor]{Proposition}
\newtheorem{propo}{Proposition}
\newtheorem{theorem}[propo]{Theorem}

 \newcommand{\cqfd}{{ \hfill
                       {\unskip\kern 6pt\penalty 500
                       \raise -2pt\hbox{\vrule\vbox to 6pt{\hrule width 6pt
                       \vfill\hrule}\vrule} \par}   }}
\title{\Large {\bf Construction and stability of a blow up solution for a nonlinear heat equation with a gradient term}}
\author{Mohamed Abderrahman Ebde\\
LAGA Universit\'e Paris 13.\\ Hatem
Zaag\\
CNRS LAGA Universit\'e Paris 13}
\date{}
\begin{document}
\maketitle
\begin{abstract}

We consider a nonlinear heat equation with a gradient term. We
construct a blow-up solution for this equation with a prescribed
blow-up profile. For that, we translate the question in selfsimilar variables and reduce the problem to a finite dimensional one. We then solve the finite dimensional problem using index theory. The interpretation of the finite dimensional parameters allows us to derive the stability of the constructed solution with
respect to initial data.
\end{abstract}
\section{Introduction}
We consider the problem
\begin{align}\label{eq:uequation}
u_{t}& =\Delta u+|u|^{p-1}u+h(u,\nabla u), \\
u(.,0)& =u(0)\in H\equiv W^{1,r}(\mathbb{R}^N)\cap W^{1,\infty}(\mathbb{R}^N),\notag
\end{align}
 $\text{ where $r$ is large enough, }  u(t):x\in \mathbb{R}^{N}\rightarrow u(x,t)\in \mathbb{R}$,
\begin{align}\label{condition:h}
       h:\mathbb R\times \mathbb {R}^N  \to \mathbb{R} &\mbox{ such that } |h(u,v)|\leq C+C|u|^{\bar{\alpha}}+ C| v|^{\alpha}, \notag\\
     p>1,~ 0\leq \bar{\alpha}<p &\mbox{ and  } 0\leq \alpha<\frac{2p}{p+1}.
     \end{align}
Note that we make no assumption on the sign of $h$.

The problem of local existence was already done by Alfonsi and
Weissler \cite{AWbugpndeap92}. They proved the uniqueness of the solution
$u(t)$ for $t\in [0,T)$, with $T\leq
\infty$. When $T<\infty$, the solution $u(t)$ blows up in the following sense:
$$\|u(t)\|_H\rightarrow +\infty \mbox{ when } t\rightarrow
T.$$
When $\alpha=\f{2p}{p+1}$, Souplet, Tayachi and Weissler \cite{STWesbugiumj96}
have managed to construct radial blow-up  selfsimilar solutions for the equation
(1),
using a shooting method to solve the differential equation.
Other  blow-up results can be found in Chipot and Weissler \cite{CWbngsj89}, Bebernes and Eberly \cite{Bebernes89}, Fila \cite{Fila91}, Kawohl and Peletier \cite{Kawohl89}.
In this paper, we focus on the subcritical case in $\alpha$ $(\alpha <\frac{2p}{p+1}),$ and aim at constructing
 a stable  blow-up  solution of the equation (\ref{eq:uequation}) and to give its behavior at  blow-up.
 Compared to the paper \cite{STWesbugiumj96}  of
Souplet, Tayachi and Weissler, our work  has a double interest:

\medskip

1-We construct a non selfsimilar solution, and give exactly its profile
which depends on the variable
$$ z=\frac{x}{\sqrt{(T-t)|\log{(T-t)}|}}.$$
Note  that the solution of
Souplet, Tayachi and Weissler, \cite{STWesbugiumj96} depends on another variable
$$y=\frac{x}{\sqrt{T-t}}.$$
2-We prove the stability of our solution with respect to perturbations in the
initial data.

\medskip

Our main result is the following:
\begin{theorem}\label{th:1}{\bf(Existence of a blow-up solution for  equation (\ref{eq:uequation}) with the
description of its profile)}
There exists $T>0$ such that equation (\ref{eq:uequation}) has a solution $u(x,t)$ in $ \mathbb{R}^N\times[0,T)$ such that:\\
i) the solution $u$ blows up in finite time T at the point $a=0$,\\
ii) \begin{equation}\label{thrm:bw-profil} \left\|(T-t)^{\frac{1}{(p-1)}}u(\cdot\sqrt{T-t},t)-f(\frac{.}{\sqrt{|\log{(T-t)}|}})\right\|_{W^{1,\infty}}\leq
\frac{C}{\sqrt{|\log(T-t)|}},
\end{equation} where
\begin{equation}\label{feq}f(z)=\left(p-1
+\frac{(p-1)^2}{4p}|z|^2\right)^{-\frac{1}{p-1}}.
\end{equation}
\end{theorem}
{\bf Remark}: We note that the profile $f$ is the same as in the nonlinear heat equation without the gradient term see
 (Bricmont and Kupiainen \cite{BKubhn94}).\\
{\bf Remark}: Note that $(i)$ follows from $(ii)$. Indeed, $(ii)$ implies that  $u(0,t) \sim \kappa (T-t)^{-\f{1}{p-1}}$  $\to +\infty$ as $t\to T,$ with \begin{equation}\label{Kappa} \kappa= f(0)=(p-1)^{-\f{1}{p-1}}.\end{equation}
{\bf Remark}: Note that estimate (\ref{thrm:bw-profil}) holds in $W^{1,\infty}(\mathbb{R})$, which a novelty of our paper. Indeed in the previous literature about neighboring equations (\cite{MZsbupdmj97}, \cite{Zihp98}, \cite{Masmoudi08}), the authors have only $L^\infty$ estimates.\\
{\bf Remark}: Note that classical methods such as energy estimates or the maximum principle break down for the equation (\ref{eq:uequation}). Indeed, there is no Lyapunov functional for (\ref{eq:uequation}), and the general form of $h(u,\nabla u)$ prevents any use of the maximum principle. To our knowledge, theorem \ref{th:1} is the first existence result of a blow-up solution for equation (\ref{eq:uequation}). As we point out in theorem {\ref{th:2}} below, this solution is stable with respect to perturbations in initial data.\\

Our approach in this paper is inspired by the method of Bricmont and Kupiainen \cite{BKubhn94},  Merle and Zaag \cite{MZsbupdmj97}
for the semilinear heat equation \begin{equation}\label{eq:semilinear} u_t= \Delta u +|u|^{p-1}u.\end{equation}
In some sense, we show in this paper that the  term $h(u,\nabla u)$ in $(\ref{eq:uequation})$ has a subcritical size when $\alpha< \f{2p}{p+1}$ and $\bar{\alpha}<p$. One may think then that our paper is just a simple perturbation of \cite{MZsbupdmj97}. If this is true in the statements, it is not the case for the proof, where we need some involved arguments to control the $h(u,\nabla u)$ term (see below Proposition \ref{prop:regu-parab-q-equation} page \pageref{prop:regu-parab-q-equation}  and Lemma \ref{lem:A_5} page \pageref{lem:A_5}). Therefore, with respect  to \cite{MZsbupdmj97}, all  that we need to do is to control the contribution $h(u,\nabla u)$. We will then refer to \cite{MZsbupdmj97} for the other contributions. However, for the reader's convenience, we will recall the main steps of the proof from \cite{MZsbupdmj97} and give details only for the new term.
We wonder whether for $\alpha=\frac{2p}{p+1}$, one can show the same result as Theorem \ref{th:1}, namely that equation (\ref{eq:uequation}) has a blow-up profile depending on the  variable $z=\frac{x}{\sqrt{(T-t)|\log (T-t)|}}$. \\
We would like to mention that the method of \cite{MZsbupdmj97}  has been successful for the following semilinear heat equation with no gradient structure:
\begin{equation}\label{eq:semilinearcomplex}
u_t= \Delta u +(1+i\delta)|u|^{p-1}u
\end{equation}
where $\delta\in {\R}$ is small and $u:{\R}^N\times[0,T)\to \C$ see \cite{Zihp98}.
Unlike our equation (1), note that  with respect to (\ref{eq:semilinear}), equation (\ref{eq:semilinearcomplex}) has the additional term  $i\delta|u|^{p-1}u$ whose size is critical with respect to the original nonlinearity $|u|^{p-1}u$. Note also that the method of \cite{MZsbupdmj97} gives the existence of a blow-up solution for the following Ginzburg-Landau equation:
\begin{equation*}
 u_t= (1+i\beta)\Delta u +(1+i\delta)|u|^{p-1}u-\gamma u
,\end{equation*} where $p>1,~ u:{\R}^N\times[0,T)\to \C$ and $ p-\delta^2-\beta\delta(p+1)>0$
(see Masmoudi and Zaag \cite{Masmoudi08}).\\
For simplicity in the parameters, we will give the proof in the particular case where $h$ is equal to its bound, namely when
\begin{equation}\label{def:h=h-min}
h(u,v)=\bar{\mu}|u|^{\bar\alpha}+\mu|v|^{\alpha}+\mu_0
\end{equation}
for some real numbers, $\bar\mu$, $\mu$ and $\mu_0$.

\medskip
The proof is done in the selfsimilar variables' framework which we introduce:
\begin{equation}\label{framwork:selfsimilar} y=\frac{x}{\sqrt{T-t}},~
  s=-\log{(T-t)},~~
  w_{T}(y,s)= (T-t)^{\frac{1}{p-1}}u(x,t))\end{equation}
where $T$ is the time  where we want the solution to blow
up.\\
For $y\in\R^N$ and $s\in [-\log T, +\infty),$   the equation of $w_{T }=w$ is
the   following
\begin{equation}\label{eq:wequation}
 \partial_s w=\Delta w-\frac{1}{2}y\cdot\nabla
  w-\dfrac{1}{p-1}w+|w|^{p-1}w
+\mu|\nabla w|^\alpha e^{-\beta s}+\bar\mu| w|^{\bar\alpha} e^{-\bar\beta s}+\mu_0e^{-\frac{p}{p-1} s}
\end{equation}
where  $\beta=\dfrac{2p-\alpha(p+1)}{2(p-1)}>0$ and $\bar\beta=\frac{p-\bar\alpha}{p-1}>0$.\\
{\bf Remark:} The fact that  $\alpha$ and $\bar\alpha $ are subcritical ( $\alpha <\frac{2p}{p+1}$ and $\bar\alpha<p$)
 is reflected in the fact that $\beta>0$ and $ \bar\beta> 0$,  which explains the little effect of the gradient term for large times. \\
\bigskip
From this transformation we see that Theorem \ref{th:1} reduces to proving the existence of a solution $w(y,s)$
for equation (\ref{eq:wequation}) such that for some $s_0\in \mathbb R$, and for all $s\geq s_0$,
$$\|w(y,s)-f(\frac{y}{\sqrt s})\|_{W^{1,\infty}}\leq \frac{C}{\sqrt s}.$$
This is "reasonable" in the sense that,
 instead of being an exact solution of (\ref{eq:wequation}), the function $f(\frac{y}{\sqrt s})$ is an approximate solution
  of (\ref{eq:wequation}) (up to the order $\frac{1}{s}$). This is clear from the fact that $f(z)$ satisfies the following

$$ \forall z\in \mathbb{R}^N,~~0=\frac{1}{2}z\cdot\nabla f(z)-\frac{f(z)}{p-1}+f(z)^{p}.$$
\bigskip
In the $w(y,s)$ framework, the proof relies on the understanding of the dynamics of the selfsimilar version of (\ref{eq:uequation})
(see equation (\ref{eq:wequation})) around the profile (\ref{feq}). We proceed in 2 steps:

\medskip
-In Step 1, we reduce the problem to a finite-dimensional problem: we show that it is enough to control a two-dimensional variable in order to control the solution (which is infinite-dimensional) near the profile.

\medskip
-In Step 2, we proceed by contradiction to solve the   finite-dimensional problem and conclude using index theory.
As in \cite{MZsbupdmj97}, \cite{Zihp98} and \cite{Masmoudi08}, it is possible to make the interpretation of the finite-dimensional variable in terms of the blow-up time and the blow-up point. This allows us to derive the stability of the profile (\ref{thrm:bw-profil}) in Theorem \ref{th:1} with respect to perturbations in the initial data. More precisely, we have the following.
\begin{theorem}\label{th:2}({\bf Stability of the solution constructed in Theorem \ref{th:1}}).
Let us denote by $\hat{u}(x,t)$ the solution constructed in Theorem \ref{th:1} and by $\hat{T}$ its blow-up time. Then,
 there exists a neighborhood $\mathcal{V}_0$ of $\hat {u}(x,0)$ in $W^{1,\infty}$ such that for any $u_0\in \mathcal{V}_0$, equation (\ref{eq:uequation}) has a unique solution $u(x,t)$ with initial data $u_0$, and $u(x,t)$ blows up in finite time $T(u_0)$ at some blow-up point $a(u_0).$ Moreover, estimate (\ref{thrm:bw-profil}) is satisfied by $u(x-a,t)$ and
\[ T(u_0)\rightarrow \hat T,~~ a(u_0)\rightarrow 0 \mbox{ as } u_0\rightarrow \hat{u}_0 \mbox{ in } W^{1,\infty}(\mathbb {R}^N)\]
\end{theorem}
The stability result follows from the reduction to a finite dimensional case as in \cite{MZsbupdmj97} for the semilinear heat equation (\ref{eq:semilinear}) with the same proof.  Since the new term $h(u,\nabla u)$ affects the  proof of  the finite dimensional  reduction and not the derivation of the stability from this reduction, we only prove the reduction and refer to \cite{MZsbupdmj97} for the stability.

\bigskip
The paper is organized as follows:\\
- In Section \ref{sec:fop}, we recall from \cite{MZsbupdmj97} the formulation of the problem.\\
- In Section \ref{sec:P-E-R}, we follow the method of \cite{MZsbupdmj97} to prove Theorem \ref{th:1}.
Note that we do not give the proof of Theorem \ref{th:2} since it follows from the proof of Theorem \ref{th:1} exactly as in the case of the semilinear heat equation \ref{eq:semilinear} treated in \cite{MZsbupdmj97}.\\
We give details only when we handle the new term $h(u,\nabla u)$.\\
\section{Formulation of the problem}\label{sec:fop}
For simplicity in the notation, we give the proof in one dimension ($N$=1), when  $h$ is given by (\ref{def:h=h-min}).
The proof remains the same for $N\geq 2$, and  for general $h$ satisfying (\ref{condition:h}).
\smallskip
We would like to find  $u_0$ initial data such that the solution $u$
of  the equation $(1)$   blows up in finite time $T$
and
$$\lim_{t\rightarrow
T}\|(T-t)^{\frac{1}{(p-1)}}u(\cdot\sqrt{T-t},t)-f(\frac{.}{\sqrt{|\log{(T-t)}|}})\|_{W^{1,\infty}}=0,\mbox{ where } f(z)=(p-1+\frac{(p-1)^2}{4p}|z|^2)^{-\frac{1}{p-1}}.$$
Using  the change of variables (\ref{framwork:selfsimilar}), this is equivalent to finding $s_0>0$ and $
w_0(y)$ such that the solution $w(y,s)$ of (\ref{eq:wequation}) with initial data
$w(y,s_0)$ satisfies
\begin{equation*} \lim_{s\rightarrow\infty}\left\|w(y,s)-f\left(\frac{y}{\sqrt s}\right)\right\|_{W^{1,\infty}}=0.
\end{equation*}
\noindent Introducing
\begin{equation}\label{eq:q+phi+def} w=\varphi+q \mbox{ where }\varphi=f\left(\frac{y}{\sqrt
s}\right)+\frac{\kappa}{2ps},\end{equation}
 where $f$ and $\kappa$  are defined in $(\ref{feq})$, and $(\ref{Kappa}),$
the problem is then reduced to constructing a function $q$ such that
\begin{equation}\label{eq:qt0} \lim_{s\rightarrow +\infty
}\|q(.,s)\|_{W^{1,\infty}}=0\end{equation}
  and $q$ is a solution of the following equation for all $(y,s)\in
\mathbb{R}\times[s_0,\infty)$
\begin{equation}\label{qequation}
\partial_s q=(\mathcal{L}+V)q+B(q)+{R}(y,s)+N(y,s),
\end{equation} where
\begin{equation}\label{eq:Operator:L}\mathcal{L}=\Delta-\frac{1}{2}y\cdot\nabla+1,~V(y,s)=p~\varphi(y,s)^{p-1}
-\frac{p}{p-1},\end{equation}
\begin{equation}\label{eq:B(q)}
B(q)=|\varphi+q|^{p-1}(\varphi+q)-\varphi^q-p\varphi^{p-1}q,\end{equation}
\begin{equation}\label{RNeq}
\begin{split}
&R(y,s)=\Delta\varphi-\frac{1}{2}y\centerdot\nabla\varphi-\frac{\varphi}{p-1}+
\varphi^p-\frac{\partial \varphi}{\partial s}\\
 \mbox{ and }& N(y,s)=\mu|\nabla\varphi+\nabla q|^\alpha e^{-\beta s}
+ \bar\mu|\varphi+ q|^{\bar\alpha}e^{-\bar\beta s}+\mu_0 e^{-\frac{p}{p-1} s}
\end{split}\end{equation}
Let $K(s,\sigma)$ be the fundamental solution of the operator $\mathcal{L}+V$.
Then, for each $ s\geq\sigma\geq~s_0$, we have:
 \begin{equation}\label{qinteg} q(s)=K(s,\sigma
)q(\sigma)+\int_{\sigma}^{s}d\tau K(s,\tau)B(q(\tau))+
 \int_\sigma ^s d\tau K(s,\tau)R(\tau)+\int_\sigma ^s d\tau K(s,\tau) N(\tau).
\end{equation}
In comparison with the case of equation $(\ref{eq:semilinear})$ treated in \cite{MZsbupdmj97}, all the terms in (\ref{qequation}) were already present in \cite{MZsbupdmj97},  except $N(y,s)$ which needs to be carefully handled. Therefore, using Lemma 3.15 page 168 and Lemma 3.16 page 169 of \cite{MZsbupdmj97}, we see that
 \begin{equation}\label{estim:B-R}
|B(q)|\leq C|q|^{\bar{p}}\mbox{ and }\|R(.,s)\|_{L^\infty}\leq \f{C}{s}\end{equation}
for $s$ large enough,
where $\bar{p}=\min(p,2).$
Using the definition (\ref{eq:q+phi+def}) of $\varphi$, we see that
\begin{equation}\label{estim:N}
|N(y,s)|\leq Ce^{-\beta s}\|\nabla q\|_{L^\infty}^\alpha+Ce^{-\bar\beta s}\| q\|_{L^\infty}^{\bar\alpha}+ Ce^{-\beta_0 s}\end{equation}
where $\beta_0=\min{(\beta,\bar{\beta})}>0$.
It is then reasonable to think that the dynamics of equation (\ref{qinteg}) are influenced by the linear part. Hence, we
first study the operator $\mathcal{L}$ (see page 543 in Bricmont and Kubiainen \cite{BKubhn94} and pages 773-775 in Abramowitz and Stegun \cite{AbramStug64}).
The operator $\mathcal{L}$    is  self-adjoint in
$D(\mathcal L) \subset L_\rho^2(\mathbb R)$ where
$$\rho(y)
=\frac{e^-\frac{|y|^2}{4}}{\sqrt{4\pi}} \mbox{ and } L_\rho^2(\mathbb{R})=\{v\in
L_{loc}^2(\mathbb R)\mbox{ such that }\int_{\mathbb R} v(y)^2\rho(y)dx<+\infty
\}.$$ \\
The spectrum of $\mathcal{L}$ is explicitly given by
$$spec(\mathcal L)=\{1-\frac{m}{2} | m\in \mathbb N\}.$$
It consists only in  eigenvalues. All the
eigenvalues are simple, and the eigenfunctions are dilation of Hermite's
polynomial: the eigenvalue $
1-\frac{m}{2}$ corresponds  to  the following eigenfunction:
\begin{equation} \label{h_m} h_m(y)=\sum_{n=0}^{[\frac{m}{2}]} \frac{m!}{n!(m-2n)!}(-1)^n
y^{m-2n}.\end{equation} Notice that  $h_m$ satisfies: \\
$\int h_n h_m \rho dx=2^n!\delta_{nm}.$
 We also introduce $k_m=h_m/\|h_m\|_{L_\rho^2(\mathbb{R})}^2$.\\
As it is mentioned in  \cite{MZsbupdmj97}, the potential  $V$ has two
fundamental properties  which will strongly  influence our strategy.\\
$(a)$ we have  $ V(.,s)\rightarrow 0$ in $L_\rho^2(\mathbb{R})$
when  $s\rightarrow +\infty.$ In practice,  the effect of $V$ in
the blow-up area $(|y|\leq C\sqrt{s})$ is regarded as a
perturbation of the effect of  $\mathcal{L}.$ \\
  $(b)$ outside of the blow-up area, we have the following property:\\
for all $\epsilon >0$, there exists $C_{\epsilon}>0$ and
$s_{\epsilon}$ such that $$ \sup_{s\geq
s_{\epsilon},~\frac{|y|}{\sqrt s}\geq
C_{\epsilon}}\left|V(y,s)-(-\frac{p}{p-1})\right|\leq \epsilon,$$
with $-\frac{p}{p-1}<-1.$ As $1$ is the largest eigenvalue of
the operator $\mathcal{L}$, outside the blow-up area  we can
consider   that  the operator $ \mathcal {L}+V$ is  an operator
with negative eigenvalues, hence, easily controlled. Considering the fact that
the behavior of $V$
is not the same inside and outside  the blow-up area, we decompose  $ q$ as
follows:

If  $\chi_0\in C_0^\infty([0,+\infty))$ with
$supp(\chi_0)\subset[0,2]$ and $\chi_0 \equiv 1 $ in $[0,1]$, we define
\begin{equation}\label{def:chi}
\chi(y,s)=\chi_0\left(\frac{|y|}{K_0~s^{\frac{1}{2}}}\right)\end{equation}
with $K_0>0$  to be fixed large enough.
If
\begin{equation}\label{def:q:proj}
q=q_b+q_e, \mbox{ with } q_b=q\chi
\mbox{ and } q_e= q(1-\chi),\end{equation} we remark that $$ \mbox{supp
}q_b(s)\subset B(0,2K_0\sqrt s),\mbox{ and supp}~ q_e(s)\subset
\mathbb{R}^N\setminus B(0,K_0\sqrt s).$$ We write
\begin{equation}\label{def:q_b} q_b(y,s)= \sum_{m=0}^2
q_m(s)h_{m}(y)+q_{-}(y,s),\end{equation}
with\\ $q_m$ is the projection of  $q_b$ on  $ h_m$ and \\
$q_{-}(y,s)= P_{-}(q_{b})$ where $P_{-}$ is the projection in the negative
subspace of the  $\mathcal{L}.$
 Thus, we can decompose $q$ in  $5$ components  as follows:\\
\begin{equation}\label{qprojection} q(y,s)=\sum_{m=0}^2
q_m(s)h_{m}(y)+q_{-}(y,s)+q_e(y,s).\end{equation}
Here and throughout the paper, we call $q_{-}
$ the negative mode of $q$, $q_0 $ the null mode of $q$, and the subspace spanned by \{$h_m / m\geq 3$\}
will be referred to as the negative subspace.
\section{Proof of the existence  result}\label{sec:P-E-R}
This section is devoted to the  proof  of the existence result (Theorem \ref{th:1}). We proceed in 4 steps, each of them making a separate subsection.

\bigskip
-In the first subsection, we define a shrinking set $V_A(s)$ and translate our goal of making $q(s)$ go to $0$ in $L^\infty(\R)$ in terms of belonging to $V_A(s)$. We also exhibit a two-parameter initial data family for equation $(\ref{qequation})$ whose coordinates are very small (with respect to the requirements of $V_A(s)$), except  the two first $ q_0$ and $q_1$.\\
-In the second subsection, we solve the local  in time Cauchy problem for equation $(\ref{qequation})$.\\
-In the third subsection, using the spectral properties  of equation $(\ref{qequation})$, we reduce our goal from the control of  $q(s)$ (an infinite-dimensional variable)  in $V_A(s)$  to the control of its two first components $(q_0,q_1)$ (a two-dimensional variable) in $[-\f{A}{s^2}, \f{A}{s^2}]^2$.\\
- In the  fourth subsection, we solve the finite-dimensional problem using index theory and conclude the proof of the theorem 1.\\
\subsection{Definition of a shrinking set $V_A(s)$ and preparation of initial data}
Let first introduce the following proposition:
\begin{prop}\label{prop:V-def}{\bf (A set shrinking to zero)} For all  $ A\geq 1$ and $ s\geq e$, we define $V_A(s)$ as
the set of all functions $r$ in
$L^\infty$ such that:\[|r_m(s)|\leq
~As^{-2},~m=0,1,~|r_2(s)|~\leq~A^2(\log s)s^{-2}
,\] \[ \forall y\in \R,~ |r_-(y,s)|~\leq~ A(1+|y^3|)s^{-2},~||r_e(s)||_{L^\infty}~\leq~
A^2s^{-\frac{1}{2}},\]
where $r_-$, $r_e$ and  $r_m$ are defined in (\ref{qprojection}). Then we have for all $s\geq e$ and $r\in V_A(s)$,

 $$ \begin{array}{lcrl}\  \mbox{(i) for all }  y\in \R,& |r_b(y,s)|&\leq& CA^2\frac{\log s}{s^2}(1+|y|^3),\\
\mbox{(ii)} &\|r(s)\|_{L^\infty}&\leq&  \frac{CA^2}{\sqrt s},\\
.\end{array}$$
\end{prop}
{\bf Proof:} (i)
For all $s\geq e$ and $r(s)\in V_A(s),$ we have
$$r_b(y,s)=\left(\sum_{m=0}^2 r_m(s)h_m(y)+r_-(y,s)\right)\cdot 1_{\{|y|\leq
2K_0\sqrt{s}\}}(y,s)$$
Using  the definition (\ref{h_m}) of $h_m$  we get:\\
\begin{equation}\label{rbestimate}
|r_b(y,s)|\leq C(1+|y|)\frac{A}{s^{2}}+ C(1+|y|^2)A^2\frac{\log s}{{s^{2}}}+ C(1+|y|^3)\frac{A}{{s^{2}}}\leq CA^2\frac{\log s}{s^2}(1+|y|^3)
 \end{equation} which gives (i).\\
(ii) Since we have \begin{align*}r(y,s)&=r_b(y,s)+ r_e(y,s),\\\end{align*}
we just use  (\ref{rbestimate}) for $(|y|\leq 2K\sqrt{s})$ together with the fact that $\|r_e\|_\infty \leq A^{2}s^{-1/2}.$ This ends the proof of Proposition \ref{prop:V-def}. \boxs\\
Initial data  (at time $s=s_0\equiv -\log T$) for the equation (\ref{qequation})
will depend on two real parameters $d_0$ and $d_1$ as given in the following proposition.

\begin{prop}\label{prop:diddc}{\bf(Decomposition of initial data on the different components)}
For each $A\geq 1$, there exists   $T_1(A)\in (0,1/e)$ such that for all $T\leq T_1$:\\
If we consider as initial data for the equation (\ref{qequation}) the following affine function :
\textit{\begin{equation}\label{qint}
 q_{d_0,d_1}(y,s_0)=f(\frac{y}{\sqrt {s_0}})^p(d_0+d_1\frac{y}{\sqrt {s_0}})-\frac{\kappa}{2ps_0}
\end{equation}}
where $f$ is defined in (\ref{feq}) and  $s_0=-\log T,$ then:\\
(i) There exists a rectangle \begin{equation}\label{D}
\mathcal{D}_T\subset [-\f{c}{s_0},\f{c}{s_0}]^2
 \end{equation}
 such that the mapping $(d_0, d_1)\rightarrow (q_0(s_0), q_1(s_0))$
(where $q$ stands $q_{d_0,d_1}(s_0)$) is linear and one to one  from $\mathcal{D}_T$ onto $[-\f{A}{s_0^2},\f{A}{s_0^2}]^2$ and maps $ \partial \mathcal{D}_T $ into $\partial [-\f{A}{s_0^2}\f{A}{s_0^2}]^2$. Moreover, it is of degree one on the boundary.
(ii)
For all $(d_0,d_1)\in \mathcal{D}_T$
we have
\begin{align}|q_2(s_0)|&\leq ~\f{\log s_0}{s_0^{2}},~
|q_{-}(y,s_0)|\leq ~\f{c}{s_0^{2}}(1+|y|^3)~ \mbox{ and }
\|q_e(.,s_0)\|_{L^\infty}\leq~\f{1}{\sqrt{s_0}},\notag\\
|d_0|+|d_1|&\leq \f{c}{s_0}\label{redinitial},\\
 \|\nabla q(\cdot,s_0)\|_{L^\infty}&\leq \frac{C_1
|d_0|}{\sqrt{s_0}}+\frac{C_2|d_1|}{\sqrt
s_0}\leq \frac{C}{\sqrt s_0} \label{qrdienq_s0}.\end{align}
$(iii)$ For all $(d_0,d_1) \in \mathcal{D}_T$,  $q_{d_0,d_1}~\in~V_A(s_0)$
with ``strict inequalities'' except for $(q_0(s_0),q_1(s_0))$, in the sense that

\[|q_m(s)|\leq
~As^{-2},~m=0,1,~|q_2(s)|<~A^2(\log s)s^{-2}
,\] \[ \forall y\in \R,~ |q_-(y,s)|<~ A(1+|y^3|)s^{-2},~||q_e(s)||_{L^\infty}<~
A^2s^{-\frac{1}{2}},\]
\end{prop}
{\bf{Proof:}}
Since we have the same definition of the set $V_A, $
and the same expression $(\ref{qint})$ of initial data $q(d_0)$ as in $\cite{MZsbupdmj97}$,
 we refer the reader to Lemma 3.5 and the Lemma 3.9  in $\cite{MZsbupdmj97}$ except for the bound $(\ref{qrdienq_s0})$
  which is new and which we prove in the following (note that although (\ref{D}), (\ref{redinitial})
   are not stated explicitly in Lemma 3.5 of \cite{MZsbupdmj97}, they are clearly written in its proof).
   In addition, for the readers convenience, we give some hint of the proof of (i).\\
   \textit{Hint of the proof of }(i)(for details, see pages 156-157 in \cite{MZsbupdmj97}): We recall form page 157 in  $\cite{MZsbupdmj97}$ the mapping
$$( d_0,d_1)\mapsto (q_0(s_0), q_1(s_0))=(d_0 a_0(s_0)+b_0(s_0), d_1 a_1(s_0)+b_1(s_0)),$$
which reduces in our case to a linear mapping :
$$( d_0,d_1)\mapsto (q_0(s_0), q_1(s_0))=(d_0 a_0(s_0), d_1 a_1(s_0)),$$
with $a_0(s_0)\neq 0$ and $a_1(s_0)\neq 0$.\\
 Therefore, it is clear that $\mathcal D_T$ is rectangle ($\mathcal D_T$ is by construction
 the set inverse image  of $[-\frac{A}{s_{0}^2}, \frac{A}{s_{0}^2}]^2$)
 and if  $(d_0, d_1) \in \partial \mathcal D_T$, then  $(q_0(s_0), q_1(s_0)) \in \partial [-\frac{A}{s_{0}^2},
  \frac{A}{s_{0}^2}]^2$.\\
\textit{Proof of $(\ref{qrdienq_s0})$}:
   Using(\ref{qint}), we have
\begin{equation*}\nabla q(y,s_0)=p\frac{d_0}{\sqrt s_0}f'(z)
f^{p-1}(z)
+\frac{d_1}{\sqrt s_0}f^{p}(z)
+p\frac{d_1}{s_0}y f'(z) f^{p-1}(z), \mbox{ where  } z=\f{y}{\sqrt s_0}
\end{equation*}
and $f$ is given in (\ref{feq}). Since $ f'(z)= -\frac{p-1}{2p}zf(z)^p,$  by (\ref{feq}), this gives,
$$ \nabla q(y,s_0)=\frac{f(z)^p}{\sqrt
s_0}(z\:p\:d_0f(z)^{p-1}+d_1+p\:d_1z^2f(z)^{p-1}).$$
 From (\ref{feq}) we can see clearly  that  $ f^p,~zf^{p-1}$ and $z^2 f^{p-1}$ are in $L^\infty(\mathbb
R)$ and we get from (\ref{redinitial}):
\begin{align*}| \nabla q(y,s_0)|\leq C_1(p)\frac{|d_0|}{\sqrt s_0}+C_2(p)\frac{|d_1|}{\sqrt s_0}\leq \frac{C}{\sqrt s_0}\end{align*} and  ($(\ref{qrdienq_s0})$) follows. This concludes the proof of Proposition \ref{prop:diddc}.
\boxs\\
Because of the presence of  $|\nabla q|^\alpha$ in the equation (\ref{qequation}), we need the following parabolic regularity estimate for equation (\ref{qequation}), with $q(s_0)$ given by (\ref{qint}) and $q(s)\in V_A(s)$. More precisely, we have the following:
\begin{prop}\label{prop:regu-parab-q-equation}{\bf(Parabolic regularity for equation (\ref{qequation}))} For all $A\geq 1$, there exists $\bar{T}_1(A)\leq T_1(A)$ such that for all $T\leq \bar{T}_1(A)$, if  $q(s)$ is a solution of equation (\ref{qequation}) on $[s_0,s_1]$ where $s_1\geq s_0=-\log T$,  $q(s_0)$ is given by (\ref{qint}) with $(d_0,d_1)\in \mathcal D_T$ and \begin{equation}\label{eq:q-in-VA-bis}
 q(s)\in V_A(s) \mbox{ for all } s\in [s_0,s_1],\end{equation}
  then, for some $C_1>0$, we have
\[ \|\nabla q(s)\|_{L^\infty}\leq \frac{C_1A^2}{\sqrt s} \mbox{ for all } s\in[s_0,s_1].\]
\end{prop}
{\bf Proof:} We consider $A\geq 1$, $T\leq T_1(A)$ and $q(s)$ a solution of equation (\ref{qequation}) defined on $[s_0,s_1]$ where $s_1\geq s_0=-\log T$ and $q(s_0)$ given by (\ref{qint}) with $(d_0,d_1)\in \mathcal D_T$. We also assume that $q(s)\in V_A(s)$ for all $s\in[s_0,s_1]$. In the following, we handle two cases: $s\leq s_0+1$ and $s\geq s_0+1.$\\
 {\bf Case 1:} $s\leq s_0+1$. Let $s_1'=\min (s_0+1,s_1)$ and take $s\in [s_0,s_1']$. Then we have for any  $t\in[s_0,s]$,
  \begin{equation}\label{canada}
 s_0\leq t\leq s\leq s_0+1\leq 2s_0, \mbox{ hence } \frac{1}{2s_0}\leq \frac{1}{s}\leq \frac{1}{t}\leq \frac{1}{s_0}
 .\end{equation}
 From equation (\ref{qequation}), we write for any   $s\in[s_0,s_1'],$
  \begin{equation}\label{eq:q:int:operator:L-case1}
  q(s)= e^{(s-s_0)\mathcal{L}} q(s_0)+ \int_{s_0}^{s} e^{(s-t)\mathcal{L}} F(t)dt,
  \end{equation}
  where \begin{equation}\label{def:F}
  F(x,t)=V(x,t)q+B(q)+R(x,t)+N(x,t)\end{equation}
  and from page 545 of \cite{BKubhn94}, for all $\theta >0 $,
  \begin{equation}\label{eq:semigroup:theta-L-case1}
   e^{\theta \mathcal{L}}(y,x)=
\frac{e^\theta}{\sqrt{4\pi(1-e^{-\theta})}}\exp{[-\frac{(ye^{-\theta/2}-x)^2}{
4(1-e^{-\theta})}]}.
  \end{equation}
   Since we easily see from (\ref{eq:semigroup:theta-L-case1}) that for any  $r\in W^{1,\infty}$ and $\theta >0$,
   \[\|\nabla (e^{\theta \mathcal{L}}r)\|_{L^{\infty}}\leq Ce^{\frac{\theta }{2}} \|\nabla r\|_{L^\infty},\]
   we write from (\ref{eq:q:int:operator:L-case1}), for all $s\in [s_0,s_1']$
   \begin{align}
   \|\nabla q(s)\|_{L^\infty}&\leq \|\nabla e^{(s-s_0)\mathcal{L}}q(s_0)\|_{L^\infty}+\int_{s_0}^{s}\|\nabla (e^{(s-t)\mathcal{L}})F(t)\|_{L^{\infty}}dt \nonumber\\
   &\leq C\|\nabla q(s_0)\|_{L^\infty}+C\int_{s_0}^{s} \frac{\|F(t)\|_{L^{\infty}}}{\sqrt{1-e^{-(s-t)}}}dt. \label{ineg:grad-q-case1}
   \end{align}
    Using (\ref{qrdienq_s0}) and (\ref{canada}), we write
    \begin{equation}\label{ineg:regular-parabo-q0}
    \|\nabla q(s_0)\|_{L^\infty}\leq \frac{C}{\sqrt {s_0}}\leq \frac{C}{\sqrt s}.\end{equation}
    Since  $\|V(t)\|_{L^{\infty}}\leq C$ (see (\ref{eq:Operator:L})), using (\ref{eq:q-in-VA-bis}), (ii) of Proposition \ref{prop:V-def}, (\ref{estim:B-R}) and (\ref{estim:N}) and (\ref{canada}), we write
   for all $t\in[s_0,s]$ and $x\in \mathbb{R}$,
  \begin{equation}\label{ineg:estim1:F}
  |F(x,t)|\leq \frac{CA^2}{\sqrt t}+e^{-\beta t}\|\nabla q(t)\|_{L^{\infty}}^{\alpha}\leq C(A^2s^{-\frac{1}{2}}+e^{-\beta s} \|\nabla q(t)\|_{L^\infty}^{\alpha}).\end{equation}
  Therefore, from (\ref{ineg:grad-q-case1}),  (\ref{ineg:regular-parabo-q0}) and (\ref{ineg:estim1:F}) we write with  $g(s)=\|\nabla q(s)\|_{L^\infty},$
  \begin{align}
  g(s)&\leq \frac{CA^2}{\sqrt s}+C\int_{s_0}^{s}\frac{ s^{-\frac{1}{2}} +e^{-\beta s} g(t)^\alpha}{\sqrt{1-e^{-(s-t)}}}dt\nonumber\\
  &\leq \frac{CA^2}{\sqrt s}+Ce^{-\beta s}\int_{s_0}^{s}\frac{g(t)^\alpha}{\sqrt{1-e^{-(s-t)}}}dt
  .\end{align}
   Using a Gronwall's argument,
  we see that for $s_0$ large enough, we have
  \[ \forall  s\in [s_0,s_1'],~ g(s)\leq \frac{2CA^2}{\sqrt s }.\]
 {\bf  Case 2: $s\geq s_0+1$}
 (note that this case does not occur when $s_1\leq s_0+1$).

 From equation (\ref{qequation}), we write for any $s'\in[s-1,s],$
  \begin{equation}\label{eq:q:int:operator:L-case2}
  q(s')= e^{(s'-s+1)\mathcal{L}} q(s-1)+ \int_{s-1}^{s'} e^{(s'-t)\mathcal{L}} F(t)dt,
  \end{equation}
where $F(x,t) $ and $e^{\theta \mathcal L}$ are given in (\ref{def:F}) and (\ref{eq:semigroup:theta-L-case1}).
Since we easily see from (\ref{eq:semigroup:theta-L-case1}) that for any  $r\in L^{\infty}$ and $\theta >0$,
   \[\|\nabla (e^{\theta \mathcal{L}}r)\|_{L^{\infty}}\leq \frac{Ce^{\frac{\theta }{2}}}{\sqrt{1-e^{-\theta}}}\|r\|_{L^\infty},\]
   we write from (\ref{eq:q:int:operator:L-case2}), for all $s'\in [s-1,s]$
   \begin{align}
   \|\nabla q(s')\|_{L^\infty}&\leq \|\nabla e^{(s'-s+1)\mathcal{L}}q(s-1)\|_{L^\infty}+\int_{s-1}^{s'}\|\nabla (e^{(s'-t)\mathcal{L}})F(t)\|_{L^{\infty}}dt \nonumber\\
   &\leq \frac{C}{\sqrt{1-e^{-(s'-s+1)}}}\|q(s-1)\|_{L^\infty}+C\int_{s-1}^{s'} \frac{\|F(t)\|_{L^{\infty}}}{\sqrt{1-e^{(s'-t)}}}dt. \label{ineg:grad-q-case2}
   \end{align}
  Since $\| q(s')\|_{L^\infty}\leq \frac{CA^2}{\sqrt {s'}}\leq \frac{CA^2}{\sqrt {s-1}}\leq \frac{CA^2}{\sqrt {s}}$ from (\ref{eq:q-in-VA-bis})
  and (ii) of Proposition \ref{prop:V-def}, we use the fact that
  $\|V(s-1)\|_{L^{\infty}}\leq C$ (see (\ref{eq:Operator:L})) and (\ref{estim:B-R}) and (\ref{estim:N}) to write for all $t\in[s-1,s']$ and $x\in \mathbb{R}$,
  \[|F(x,t)|\leq \frac{CA^2}{\sqrt t}+e^{-\beta t}\|\nabla q(t)\|_{L^{\infty}}^{\alpha}\leq C((A^2{s'}^{-\frac{1}{2}}+e^{-\beta s'} \|\nabla q(t)\|_{L^\infty}^{\alpha}).\]
  Therefore, from (\ref{ineg:grad-q-case2}), we write with  $g(s')=\|\nabla q(s')\|_{L^\infty}$
  \begin{align}
  g(s')&\leq \frac{CA^2}{\sqrt {s'} \sqrt {1-e^{-(s'-s+1)}}}+C\int_{s-1}^{s'}\frac{ {s'}^{-\frac{1}{2}} +e^{-\beta {s'}} g(t)^\alpha}{\sqrt{1-e^{-(s'-t)}}}dt\nonumber\\
  &\leq \frac{CA^2}{\sqrt {s'}\sqrt {1-e^{-(s'-s+1)}}}+Ce^{-\beta {s'}}\int_{s-1}^{s'}\frac{g(t)^\alpha}{\sqrt{1-e^{-(s'-t)}}}dt
  .\end{align}
   Using a Gronwall's argument,
  we see that for $s$ large enough,
  \[ \forall  s'\in [s-1,s],~ g(s')\leq \frac{2CA^2}{\sqrt {s'} \sqrt {1-e^{-(s'-s+1)}}}.\]
 Taking $s'=s$ concludes the proof of Proposition \ref{prop:regu-parab-q-equation}.\Box

\subsection[]{Local in time solution of  equation (\ref{qequation})}
In the following, we find a local in time solution for equation (\ref{qequation}).
\bigskip
\begin{prop}\label{prop:localq}{\bf(Local in time solution for equation (\ref{qequation}) with initial data (\ref{qint}))}
For all $A\geq 1$, there exists $T_2(A)\leq \bar{T}_1(A)$ such that for all $T\leq T_2$, the following holds:\\
For all $(d_0, d_1)\in\mathcal{D}_T $ with initial data at $s=s_0$ $q_0(s_0) $ defined in (\ref{qint}), there exists
 a $s_{max}>s_0$ such that  equation (\ref{qequation}) has a unique solution satisfying $q(s)\in V_{A+1}(s)$ for all $s\in [s_0, s_{max})$.
\end{prop}
{\bf Proof.} From the definition of $q$ in (\ref{eq:q+phi+def}) we can see that the Cauchy problem of (\ref{qequation}) is equivalent to that of equation (\ref{eq:wequation}) which is equivalent to the Cauchy problem of equation (\ref{eq:uequation}).
 Moreover, the initial data $q_0$ defined in (\ref{qint}) gives  the following  initial data of (\ref{eq:uequation}):
$$ u_{0,d_1,d_2}(x)=T^{-\frac{1}{(p-1)}}\left\{f(z)\left(1+\frac{d_0+d_1z}{
p-1+\frac{(p-1)^2}{4p}z^2}\right)\right\}, \mbox{ where } z=x(|\log{T}|T)^{-\frac{1}{2}}.$$
  which  belongs to $H\equiv W^{1,r}(\mathbb{R}^N)\cap W^{1,\infty
}(\mathbb{R}^N)$ for $r$  which insures the local existence (see  the introduction) of $u$ in $H$.
Now, since we have from (iii) of Proposition \ref{prop:diddc}, $q_0\in V_A(s_0)\subsetneqq V_{A+1}(s_0)$, there exists $s_3$ such that for all $s\in [s_0, s_3)$, $q(s)\in V_{A+1}(s)$. This concludes the proof of Proposition \ref{prop:localq}.
\boxs
\bigskip
\subsection{ Reduction to a finite-dimensional problem}
\bigskip
 This step is crucial in the proof of Theorem \ref{th:1}. In this step,
we will prove through a priori estimates that for each
 $s\geq s_0$, the control of
$q(s)$ in   $V_A(s)$  is reduced to the control of $ (q_0,q_1)(s)$ in
 $ [-\f{A}{s^2},\f{A}{s^2}]^2$.
In fact,  this result implies that if for some
$s_1\geq s_0,~~q(s_1)\in\partial V_A(s_1)$ then
$(q_0(s_1),q_1(s_1))\in \partial[-\f{A}{s^2},\f{A}{s^2}]^2.$\\

\begin{prop}\label{prop:rt2dimen} {\bf(\noindent Control of  $q(s)$ by $(q_0(s),q_1(s) )$ in $ V_A(s) $)} There
exists $A_3>0$ such that  for each $A\geq~A_3$
there exists $T_3(A)\leq T_2(A) \mbox{ such that for all }T\leq T_3(A)$,
 the following holds: \\
If $q$ is a solution of (\ref{qequation}) with initial data at $s=s_0=-\log T$ given by
(\ref{qint}) with $(d_0,d_1)~\in~\mathcal{D}_T $, and $q(s)~\in~V_A(s)$  for all  $ s \in~[s_0,s_1],~$
with $q(s_1)\in \partial V_A(s_1)$ for some $s_1\geq s_0$, then:\\
(i) (Reduction to a finite dimensional problem) $ (q_0(s_1),q_1(s_1))\in \partial [-\f{A}{s_1^2},\f{A}{s_1^2}]^2$.\\
(ii) (Transverse crossing) There exist $m\in \{0,1\}$ and $\omega \in \{-1,1\}$ such that
\begin{equation*}
 \omega q_m(s_1)=\f{A}{s_1^2}\mbox{ and }  \omega\f{d q}{d s}(s_1)>0.
\end{equation*}
\end{prop}
{\bf Remark.} In $N$  dimensions, $q_0\in \R$ and $q_1\in \R^N$. In particular, the finite-dimensional problem is of dimension  $N+1$. This is why in initial data (\ref{qint}), one has to take $d_0\in \R$ and $d_1\in \R^N$.\\
{\bf Proof:}
Let us consider $A\geq 1$ and $T\leq T_2(A)$. We then consider $q$ a solution of $(\ref{qequation})$ with initial data at $s=s_0=-\log T$ given by (\ref{qint}) with $(d_0,d_1)\in\mathcal D_T$,  and $q(s)\in V_A(s)$ for all $s\in[s_0,s_1]$, with $q(s_1)\in \partial V_A(s_1)$ for some  $s_1\geq s_0$.
 We now claim the following:

\begin{prop}\label{prop:projet:q}
There exists $A_4\geq 1$ such that for all $A\geq A_4$ and $\rho\geq 0,$ there exists $s_4(A,\rho)\geq -\log T_2(A)$ such that the following holds for all $s_0\geq s_4(A,\rho)$:\\
Assume that for all $s\in [\tau,\tau+\rho]$ for some
$\tau\geq s_0,$
\[q(s)\in V_A(s)\mbox{ and } \|\nabla q(s)\|_{L^\infty}\leq \frac{CA^2}{\sqrt s}.\] Then, the following holds for all $s\in[\tau,\tau+\rho]$:\\
$(i)$ (ODE satisfied by the expanding modes) For $m=0$ and $1$, we have
$$ \left|q_m'(s) -(1-\frac{m}{2})q_m(s)\right|\leq \frac{C}{s^2}.$$
$(ii)$ (Control of the null and negative modes) we have \\
 \begin{align*}
|q_2(s)|&\leq \frac{\tau^2}{s^2}|q_{2}(\tau)|+\frac{C A(s-\tau)}{s^3},\\
\left\| \frac{q_{-}(s)}{1+|y|^3}\right\|_{L^{\infty}}&\leq C e^{-\frac{(s-\tau)}{2}}\left\|\frac{q_{-}(\tau)}{1+|y|^3}\right\|_{L^\infty}
+C\frac{e^{-(s-\tau)^2}\| q_{e}(\tau)\|_{L^\infty}}{s^{3/2}}+\frac{C(s-\tau)}{s^2},\\
\|q_{e}(s)\|_{L^\infty}&\leq e^{-\frac{(s-\tau)}{p}}\|q_{e}(\tau)\|_{L^\infty}+e^{s-\tau} s^{3/2}\left\|\frac{q_{-}(\tau)}{1+|y|^3}\right\|_{L^\infty}+\frac{(s-\tau)}{s^{1/2}}.
\end{align*}

\end{prop}
 {\bf Proof:} See Section \ref{sec:pfdp}.
\boxs\\
 Using Propositions \ref{prop:projet:q} and \ref{prop:diddc},  one can see that Proposition \ref{prop:rt2dimen} follows exactly as in the case of semilinear equation (\ref{eq:semilinear}) treated in \cite{MZsbupdmj97}. The proof is easy, however, a bit technical. That is the reason why it is  omitted. The interested reader can find details in pages (164-160) of \cite{MZsbupdmj97} for $(i)$, and in page 158 of \cite{MZsbupdmj97} for $(ii)$.
 This concludes the proof of Proposition \ref{prop:rt2dimen}.
 \boxs
\subsection{Proof of the finite dimensional problem and proof of Theorem \ref{th:1}}
We prove Theorem \ref{th:1} using the previous subsections. We proceed in  two parts.\\

-In Part 1, we solve the finite-dimensional problem and prove the existence of $A\geq 1,~ T>0,$ $(d_0,d_1)\in \mathcal{D}_T$ such that problem (\ref{qequation}) with initial data at $s=s_0=-\log T,~q_{d_0,d_1}(s_0)$ given by
(\ref{eq:qt0}) has a solution $q(s)$ defined for all $s\in[-\log T,\, \infty)$ such that
\begin{equation}\label{eq:q-in-VA}
q(s)\in V_A(s) \mbox{ for all } s\in[-\log T, \,+\infty).
\end{equation}

-In Part 2, we show that the solution constructed in Part 1 provides a blow-up solution of equation (\ref{eq:uequation})
which satisfies all the properties stated in Theorem \ref{th:1}, which concludes the proof.\\
{\bf Part 1: Solution  of the 2-dimensional problem}\\
This part in not new. Indeed, given Proposition \ref{prop:rt2dimen}, the solution of the 2-dimensional problem is exactly the same as in \cite{MZsbupdmj97} (see also \cite{Masmoudi08}). Nevertheless, for the sake of completeness we give a sketch of the argument here.
We take $A=A_3$ and $T=T_3(A)=\min(T_1(A),\bar{T}_1(A),  T_2(A), T_3(A))$ so that Propositions \ref{prop:diddc}, \ref{prop:regu-parab-q-equation}, \ref{prop:localq} and \ref{prop:rt2dimen} apply. We will find the parameter $(d_0,d_1)$ in the set
$\mathcal{D}_T$ defined in Proposition \ref{prop:diddc}. We proceed by contradiction and assume from
(iii) of Proposition \ref{prop:diddc} that for all $(d_0,d_1)\in \mathcal{D}_T$, there exists $s_*(d_0,d_1)\geq -\log T$ such that $q_{d_0,d_1}(s)\in V_A(s)$ for all $s\in[-\log T, s_*]$ and $q_{d_0,d_1}(s_*)\in \partial V_A(s_*)$.
From (ii) of Proposition \ref{prop:rt2dimen}, we see that $(q_0(s_*), q_1(s_*))\in \partial[-\frac{A}{s_*^2},\frac{A}{s_*^2}]^2$ and the following function is well defined:
 \begin{align*}\Phi(y,s)&: \mathcal{D}_{T} \rightarrow \partial [-1,1]^2\\
&(d_0,d_1)\mapsto
\dfrac{s_*^2}{A}(q_0,q_1)_{(d_0,d_1)}(s_*).\end{align*}
From (ii) of Proposition \ref{prop:rt2dimen}, $\Phi$ is continuous. If we manage to prove that $\Phi$ is of degree 1 on the boundary, then we have a contradiction from the degree theory.
Let us prove that.\\
Using (i) and (iii) of Proposition \ref{prop:diddc}  and the fact that $q(-\log T)= q_{d_0,d_1}$, we see that when $(d_0,d_1)$ is on the boundary of the rectangle $\mathcal{D}_T$, $(q_{0}, q_1)(-\log T)\in\partial[-\frac{A}{(\log T)^2},\frac{A}{(\log T)^2}]^2$ and $q(-\log T)\in V_A(-\log T)$ with strict inequalities for the other components.
Applying the transverse crossing property of (ii) in Proposition \ref{prop:rt2dimen}, we see that
$q(s)$  leaves  $V_A(s)$ at $s=-\log T$,  hence $s_*=-\log T$. Using (i) of Proposition \ref{prop:diddc}, we see that  restriction of $\phi$
to the boundary is of degree 1. A contradiction then follows. Thus, there exists a value $(d_0,d_1)\in \mathcal{D}_T$
such that for all $s\geq -\log T,~ q_{d_0, d_1}(s)\in V_A(s)$.
This concludes the proof of (\ref{eq:q-in-VA}).\\

 {\bf Part 2: Proof of Theorem \ref{th:1}:}\\
 Since $q_{d_0,d_1}$ satisfies (\ref{eq:q-in-VA}), we clearly see, from section \ref{sec:fop} and (ii) in Proposition \ref{prop:V-def} that $u$ defined from $q_{d_0,d_1}$ through the transformations (\ref{eq:wequation}) and  (\ref{eq:q+phi+def}) is well defined for all  $(x,t)\in \mathbb{R}\times[0, T)$ and satisfies (\ref{thrm:bw-profil}) in the $L^\infty$ norm.
 Using the parabolic estimate of Proposition \ref{prop:regu-parab-q-equation}, we see that (\ref{thrm:bw-profil}) holds in the $W^{1,\infty}$ norm as well. In particular from (\ref{thrm:bw-profil}), we have  $$u(0,t)\sim \kappa (T-t)^{-\frac{1}{p-1}},$$ hence, $u$ blows up at time $T$ at the point $a=0$. This concludes the proof of Theorem\ref{th:1}.\Box
 \section{Proof Proposition \ref{prop:projet:q}}\label{sec:pfdp}
In this section, we prove  Proposition \ref{prop:projet:q}.
In the following, we will  find the main contribution in the projection given in the decomposition (\ref{qprojection})
 of the four terms appearing in the right-hand side of equation (\ref{qequation}).\\

Let us recall the equation (\ref{qinteg}) of $q$\\
\begin{equation}q(s)=K(s,\tau )q(\tau)+\int_{\tau}^{s}d\sigma
K(s,\sigma)B(q(\sigma))+
 \int_\tau ^s d\sigma K(s,\sigma)R(\sigma)+\int_\tau ^s d\sigma K(s,\sigma)
N(\sigma)\end{equation}
where $K$ is the fundamental solution  of the operator  $\mathcal{L }+V.$\\ We
write
   $ q=\alpha+\beta+\gamma+\delta $ where
\begin{align} \alpha(s)&=K(s,\tau )q(\tau),~ \beta(s)=\int_{\tau}^{s}d\sigma
K(s,\sigma)B(q(\sigma)),\\
\gamma(s)&=\int_\tau ^s d\sigma K(s,\sigma)R(\sigma), \mbox{ and }
\delta(s)=\int_\tau ^s d\sigma K(s,\sigma) N(\sigma).\label{eq:delta:eq} \end{align}
 We assume that
$q(s)$ is in
$V_A(s)$ for each
 $s\in [\tau,\tau+\rho]$. Using $(41)$, we derive new bounds on
$\alpha,\,\beta,\,\gamma \mbox{ and }\tau.$

Clearly, Proposition \ref{prop:projet:q} follows from the following:
\begin{lem}\label{lem:A_5} There exists $A_5>0$ such that for all $A\geq A_5
,$ and $\rho>0$
there exists $T_5(A,\rho)\leq T_2(A)$, such that
for all $T\leq T_5(A,\rho),$  if we
assume that
for all $s\in[\tau, \tau+\rho]$, $ q(s)$ satisfies (\ref{qequation}), $q(s)\in V_A(s)$ and
 $\|\nabla q(s)\|_{L^\infty}\leq \frac{CA^2}{\sqrt s}$ with $\tau\geq s_0=-\log T.$
Then, we have the following results for all $s\in[\tau, \tau+\rho]$:\\
$(i)$\noindent{\bf(Linear term)}\\
\begin{align*}|\alpha_2(y,s)|&\leq \frac{\tau^2}{s^2}|q_2(\tau)|+(s-\tau)CAs^{-3},\\
\left\|\frac{\alpha_{-}(y,s)}{1+|y|^3}\right\|_{L^{\infty}}&\leq Ce^{-\frac{1}{2}(s-\tau)}\left\|\frac{q_{-}(\tau)}{1+|y|^3} \right\|_{L^{\infty}}+Ce^{-(s-\tau)^2}
\frac{\|q_e(\tau)\|_{L^\infty}}{s^{3/2}},\\
\|\alpha_e(s)\|_{L^\infty}&\leq C e^{-\frac{(s-\tau)}{p}}\|q_e(\tau)\|_{L^\infty}+Ce^{s-\tau}\left\|\frac{q_{-}(\tau)}{1+|y|^3} \right\|_{L^{\infty}}s^
{-\frac{3}{2}}.
\end{align*}
$(ii)$ \noindent{\bf (Nonlinear term)}\\
\begin{equation*}|\beta_2(s)|~\leq\frac{(s-\tau)}{s^{3}},~
|\beta_-(y,s)|\leq~(s-\tau)(1+|y|^3)s^{-2},~
\|\beta_{e}(s)\|_{L^\infty}~\leq~(s-\tau)s^{-\frac{1}{2}},\end{equation*}
$(iii)$  \noindent{\bf (Corrective term)}\\
$$
\begin{array}{rclrcl} |\gamma_2(s)| & \leq & C(s-\tau)s^{-3}, &
|\gamma_-(y,s)| & \leq& C(s-\tau)(1+|y|^3)s^{-2},\\
\|\gamma_e(s)\|_{L^\infty}&\leq& (s-\tau)s^{-1/2}, & & &\end{array} $$
$(iv)$ \noindent{\bf(New perturbation term)}
$$\begin{array}{lcrlcrlcr}|\delta_2(s)| &\leq& C(s-\tau)s^{-3},&
|\delta_{-}(y,s)|&\leq&C(s-\tau)(1+|y|^3)s^{-3},&
\|\delta_e(s)\|_{L^\infty}&\leq&C(s-\tau)s^{-3}\end{array}.$$\end{lem}
 \emph{ {Proof of Lemma \ref{lem:A_5} }}:
 We consider $A\geq 1$, $\rho>0$ and $T\leq e^{-\rho}$ (so that  $s_0=-\log T\geq \rho$). We then consider $q(s)$ a solution of (\ref{qequation}) satisfying
 \begin{equation}\label{q-in-V_A-and-nabla-q-samll}
 q(s)\in V_A(s) \mbox{ and } \|\nabla q(s)\|_{L^\infty}\leq \frac{CA^2}{\sqrt s } \mbox{ for all } s\in[\tau,\tau+\rho]
 \end{equation}
  for some $\tau\geq s_0$.
  The terms $\alpha$, $\beta$ and $\gamma$ are already present with case of the semilinear heat equation (\ref{eq:semilinear}), so we  refer to Lemma 3.13 page 167 in \cite{MZsbupdmj97} for their proof, and we only focus on the new term $\delta(y,s)$.
  Note that since  $s_0\geq \rho$, if  we take $\tau\geq s_0$, then $\tau+\rho\leq 2\tau$ and if  $\tau\leq \sigma\leq s\leq \tau+\rho$ then
  \begin{equation}\label{order:inverse:tau-t:s}
  \frac{1}{2\tau}\leq\frac{1}{s}\leq\frac{1}{\sigma}\leq \frac{1}{\tau}.\end{equation}
 Using (\ref{estim:N}), (\ref{q-in-V_A-and-nabla-q-samll}), and (ii) of Proposition {\ref{prop:V-def}}, we write for
 $\tau\leq \sigma \leq s\leq \tau+\rho$ and $x\in \mathbb R$
 \begin{equation}\label{estim:new-estimate:N}
 |N(x,t)|\leq Ce^{-\beta \sigma}\frac{A^{2\alpha}}{\sigma^{\frac{\alpha}{2}}}+Ce^{-\bar{\beta}\sigma}\frac{A^{2\alpha}}{\sigma^{\frac{\alpha}{2}}}+e^{-\beta_0 \sigma}\leq \frac{1}{\sigma^4}\leq \frac{C}{s^4}\end{equation}
 for $s_0$ large enough (that is $T$ small enough), where we used (\ref{estim:N}) in the last inequality.
 Recalling from Bricmont and Kupiainen \cite{BKubhn94} that for all $y,x \in \mathbb{R}$
\begin{equation}\label{ineq:B-K} |K(s,\sigma,y,x)|~\leq Ce^{(s-\sigma)\mathcal{L}}(y,x)\end{equation}
 where $e^{\theta \mathcal{L}}$ is given in (\ref{eq:semigroup:theta-L-case1})
(see page 545 in \cite{BKubhn94}), we write from
(\ref{eq:delta:eq}), (\ref{ineq:B-K}), \ref{estim:new-estimate:N}  and
(\ref{eq:semigroup:theta-L-case1}),
\begin{align*}|\delta (y,s)|&\leq \int_{\tau}^{s} d\sigma\int_{\mathbb{R}} K(s,\sigma,y,x)|N(x,\sigma)|dx\\
&\leq \int_{\tau}^{s} d\sigma \int_{\mathbb{R}} e^{(s-\sigma)\mathcal L}(y,x)\frac{C}{s^4}dx
\leq   (s-\tau)e^{(s-\tau)}\frac{C}{s^4}\\
& \leq C(s-\tau)\frac{e^\rho}{s^4}\leq\frac{(s-\tau)}{s^3}\end{align*}
for $s_0$ large enough.\\
By definitions (\ref{def:q:proj}) and (\ref{def:q_b})
 of $q_m$, $q_-$ and $q_e$, we write  for $m\leq 2$,
 \begin{equation*}\label{estim:delta_m}
 |\delta_m(s)|\leq \left|\int_{\mathbb R}|\chi(y,s)\delta(y,s)k_m(y)\rho(y)dy\right|\leq C\int_{\mathbb{R}}|\delta(y,s)|(1+|y|^2)\rho(y)dy
 \leq \frac{C}{s^3}(s-\tau)\end{equation*}
 \begin{equation*}\label{estim:delta_-}
 |\delta_-(y,s)|=|\delta(y,s)-\sum_{i=0}^2\delta_i(y,s)k_i(y)|\leq  (s-\tau)(1+|y|^3)\frac{C}{s^3}
 \end{equation*}
 \begin{equation*}\label{estim:delta_e}
  |\delta_e(y,s)|=|(1-\chi(y,s))\delta(y,s)|\leq |\delta(y,s)|\leq (s-\tau)\frac{C}{s^3},
 \end{equation*}
which concludes the proof of Lemma \ref{lem:A_5}  and Proposition \ref{prop:projet:q} too.

\begin{flushright}$\blacksquare$\end{flushright}

{\bf Acknowledgments:}{
 This work was supported by the Marie Curie Actions of the
European Commission in the frame of the DEASE project (MEST-CT-2005-021122).

The first author wishes to thank Prof. Dr. Christian Schmeiser for his hospitality in WPI in Vienna  where a part of this work was done.}
\bibliographystyle{plain}

\
\noindent  Mohamed Abderrahman EBDE:\\ Universit\'e Paris 13, Institut Galil\'ee,
Laboratoire Analyse, G\'eom\'etrie et Applications, CNRS UMR 7539, 99
avenue J.B. Cl\'ement, 93430 Villetaneuse, France.\\
e-mail: ebde@math.univ-paris13.fr

 \noindent Hatem Zaag:\\ Universit\'e Paris 13, Institut Galil\'ee,
Laboratoire Analyse, G\'eom\'etrie et Applications, CNRS UMR 7539, 99
avenue J.B. Cl\'ement, 93430 Villetaneuse, France. \\
e-mail: Hatem.Zaag@univ-paris13.fr

\end{document}